\documentclass{article}
\usepackage{color}
\usepackage{amsmath,amssymb}
\catcode`\@=11 \@addtoreset{equation}{section}

\catcode`\@=12

\newtheorem{Theorem}{Theorem}[section]
\newtheorem{Proposition}{Proposition}[section]
\newtheorem{Lemma}{Lemma}[section]
\newtheorem{Corollary}{Corollary}[section]
\newtheorem{Definition}{Definition}[section]
\newtheorem{Remark}{Remark}[section]

\newcommand{\bTheorem}[1]{
\begin{Theorem} \label{T#1} }
\newcommand{\eT}{\end{Theorem}}

\newcommand{\bProposition}[1]{
\begin{Proposition} \label{P#1}}
\newcommand{\eP}{\end{Proposition}}

\newcommand{\bLemma}[1]{
\begin{Lemma} \label{L#1} }
\newcommand{\eL}{\end{Lemma}}

\newcommand{\bCorollary}[1]{
\begin{Corollary} \label{C#1} }
\newcommand{\eC}{\end{Corollary}}

\newcommand{\bFormula}[1]{
\begin{equation} \label{#1}}
\newcommand{\eF}{\end{equation}}

\newcommand{\Om}{\Omega}

\newcommand{\DC}{C^\infty_c}

\newcommand{\vr}{\varrho}

\newcommand{\vrt}{\tilde{\vr}}

\newcommand{\vut}{\tilde{\vu}}

\newcommand{\vu}{\vc{u}}
\newcommand{\vz}{\vc{z}}
\newcommand{\vU}{\vc{U}}
\newcommand{\vc}[1]{{\bf #1}}
\newcommand{\qed}{\bigskip \rightline {Q.E.D.} \bigskip}
\newcommand{\Div}{{\rm div}_x}
\newcommand{\Grad}{\nabla_x}

\newcommand{\tn}[1]{\mbox {\F #1}}
\newcommand{\dx}{{\rm d} {x}}
\newcommand{\vph}{{\boldsymbol \varphi}}

\newcommand{\dt}{{\rm d} t }

\newcommand{\dxdt}{\dx\dt}

\newcommand{\bProof}{{\bf Proof: }}

\newcommand{\vV}{\vc{V}}
\newcommand{\n}{\vc{n}}

\newcommand{\Ecal}{{\cal E}}
\newcommand{\Rcal}{{\cal R}}
\newcommand{\ep}{\varepsilon}

\font\F=msbm10 scaled 1000

\providecommand{\abs}[1]{\left\lvert#1\right\rvert}
\newcommand{\sil}{\rightarrow}
\newcommand{\norm}[1]{\ensuremath{\left\| #1 \right\|}}
\newcommand{\norma}[2]{\ensuremath{\norm{#1}_{#2}}}

\definecolor{grey}{rgb}{0.85,0.85,0.85}
\date{}

\long\def\greybox#1{%
    \newbox\contentbox%
    \newbox\bkgdbox%
    \setbox\contentbox\hbox to \hsize{%
        \vtop{
            \kern\columnsep
            \hbox to \hsize{%
                \kern\columnsep%
                \advance\hsize by -2\columnsep%
                \setlength{\textwidth}{\hsize}%
                \vbox{
                    \parskip=\baselineskip
                    \parindent=0bp
                    #1
                }%
                \kern\columnsep%
            }%
            \kern\columnsep%
        }%
    }%
    \setbox\bkgdbox\vbox{
        \color{grey}
        \hrule width  \wd\contentbox %
               height \ht\contentbox %
               depth  \dp\contentbox
        \color{black}
    }%
    \wd\bkgdbox=0bp%
    \vbox{\hbox to \hsize{\box\bkgdbox\box\contentbox}}%
    \vskip\baselineskip%
}

\begin{document}

%%%%%%%%%%%%%%%%%%%%%%%%%%%%%%%%

\title{Weak-strong uniqueness for compressible Navier-Stokes system with slip boundary conditions on time dependent domains}

\author{
Ond{\v r}ej Kreml $^{1}$\thanks{The research of O.K. and \v S.N. is part of the MOBILITY project 7AMB16PL060 in the general framework of RVO:67985840.}
\and {\v S}{\' a}rka Ne{\v c}asov{\' a} $^1$\footnotemark[1]
\and Tomasz Piasecki $^2$\thanks{The research of T.P. was supported by Polish NCN grant UMO-2014/14/M/ST1/00108.}
}

\maketitle

\bigskip

\centerline{$^1$Institute of Mathematics of the Academy of Sciences of
the Czech Republic} \centerline{\v Zitn\' a 25, 115 67 Praha 1,
Czech Republic}
\bigskip

\centerline{$^2$  Institute of Mathematics, Polish Academy of Sciences}
\centerline{\'Sniadeckich 8, 00-656 Warszawa, Poland}

\medskip

\begin{abstract}

We consider the compressible Navier-Stokes system on 
time-dependent domains with prescribed motion of the boundary, 
supplemented with slip boundary conditions for the velocity. We derive the
relative entropy inequality in the spirit of \cite{FeJiNo} for the system on moving domain 
and use it to prove the weak-strong uniqueness property.

\end{abstract}

\medskip


\medskip

\section{Introduction}\label{i}

The flow of a compressible barotropic viscous fluid is in the absence of external forces described by the following system of partial differential equations
\bFormula{i1a}
\partial_t \vr + \Div (\vr \vu) = 0,
\eF
\bFormula{i1b}
\partial_t (\vr \vu) + \Div (\vr \vu \otimes \vu) + \Grad p(\vr) = \Div \tn{S}(\Grad \vu).
\eF
Here $\vr$ denotes the density of the fluid and $\vu$ is the velocity. The stress tensor $\tn{S}$ is determined by the standard Newton rheological law
\bFormula{i4}
\tn{S} (\Grad \vu) = \mu \left( \Grad \vu +
\Grad^t \vu - \frac{2}{3} \Div \vu \tn{I} \right) + \eta \Div \vu \tn{I}
\eF
with $\mu > 0$ and $\eta \geq 0$ being constants.

The pressure $p(\vr)$ is a given continuously differentiable function of the density. We introduce the pressure potential as
\bFormula{p1pb}
H(\vr) = \vr\int_1^\vr \frac{p(z)}{z^2} {\rm d} z.
\eF

We are interested in properties of solutions to the system of equations \eqref{i1a}-\eqref{i1b} on a moving domain $\Omega = \Omega_t$ with prescribed movement of the boundary. More precisely, the boundary of the domain $\Omega_t$ occupied by the fluid is assumed to be described by means of a given velocity field $\vc{V}(t,x)$, where $t \geq 0$ and
$x \in \mathbb{R}^3$. Assuming $\vc{V}$ is regular, we solve the associated system of differential equations
\begin{equation}\label{eq:coc}
\frac{{\rm d}}{{\rm d}t} \vc{X}(t, x) = \vc{V} \Big( t, \vc{X}(t, x) \Big),\ t > 0,\ \vc{X}(0, x) = x,
\end{equation}
and set
\[
\Omega_\tau = \vc{X} \left( \tau, \Omega_0 \right), \ \mbox{where} \ \Omega_0 \subset \mathbb{R}^3 \ \mbox{is a given domain},\
\Gamma_\tau = \partial \Omega_\tau, \ \mbox{and}\ Q_\tau = \{ (t,x) \ |\ t \in (0,\tau), \ x \in \Omega_\tau \}.
\]

The impermeability of the boundary of the physical domain is described by the condition
\bFormula{i6} (\vu - \vc{V})
\cdot \vc{n} |_{\Gamma_\tau} = 0 \ \mbox{for any}\ \tau \geq 0,
\eF
where $\vc{n}(t,x)$ denotes the unit outer normal vector to the boundary $\Gamma_t$.

Moreover, we prescribe the Navier slip boundary conditions, i.e.
\bFormula{b1} 
\left(\left[ \tn{S} \vc{n} \right]_{\rm tan} + \kappa\left[ \vu - \vc{V} \right]_{\rm tan}\right)|_{\Gamma_\tau} = 0,
\eF
where $\kappa \geq 0$ represents a friction coefficient. In particular the choice $\kappa = 0$ yields the complete slip boundary condition.
%
%For physical motivation of correct description of the fluid
%boundary behavior, see Bul\' \i \v cek, M\' alek and Rajagopal
%\cite{BMR1}, Priezjev and Troian \cite{PRTR} and the references
%therein.

Finally, the system of equations \eqref{i1a}-\eqref{i1b} is supplemented by the initial conditions
\bFormula{i7} 
\vr(0, \cdot) = \vr_0 ,\quad 
(\vr \vu) (0, \cdot) = (\vr \vu)_0 \quad \mbox{in}\ \Omega_0. 
\eF

Global weak solutions to the compressible barotropic Navier--Stokes system on a fixed domain were proved to exist in a pioneering work by Lions \cite{LI4}. This theory was later extended by Feireisl and collaborators (\cite{FNP}, \cite{EF61}, \cite{EF70}, \cite{EF71}) to cover larger class of pressure laws. The existence theory in the case of moving domains was developed in \cite{FeNeSt} for no-slip boundary conditions using the so-called Brinkman penalization and in \cite{FKNNS} for slip boundary conditions.

There is also quite large amount of literature available concerning existence of strong solutions for the system \eqref{i1a}--\eqref{i1b} (or even for more complex systems involving also heat conductivity of the fluid) on a fixed domain, such solutions are proved to exist either locally in time or globally provided the initial data are sufficiently close to a rest state, let us mention among others \cite{MN1}, \cite{MN2}, \cite{V1}, \cite{V2}. All these results are proved under the no-slip boundary conditions. 
The only global existence result without smallness assumptions on the data has been obtained by Vaigant and Kazhikov \cite{VK} 
for space periodic boundary conditions in two space dimensions.
The case of slip boundary condition on fixed domain was considered by Zajaczkowski \cite{Za} and Hoff \cite{Ho}. Even local-in-time existence results of strong solutions on moving domains (regardless of boundary behavior) seems to be an open question at this moment and it will be a matter of our forthcoming paper \cite{KNP2}.

The concept of relative entropies has been succesfully used in the context of partial differential equations (see among others Carillo et al. \cite{Ca}, Masmoudi \cite{MA}, Saint-Raymond \cite{SR}, Wang and Jiang \cite{WaJa}). Germain \cite{Ge} introduced a notion of solution to the system \eqref{i1a}--\eqref{i1b} based on a relative entropy inequality with respect to a hypothetical strong solution. Similar idea was adapted by Feireisl et al. \cite{FeNoSu} who defined a suitable weak solution to the barotropic Navier-Stokes system based on a general relative entropy inequality with respect to any sufficiently smooth pair of functions. In \cite{FeJiNo} the authors used relative entropy inequality to prove the weak-strong uniqueness property. Doboszczak \cite{dobo} proved both the relative entropy inequality as well as the weak-strong uniqueness property in the case of moving domain and no-slip boundary condition.

In this paper we first prove that weak solutions to the system \eqref{i1a}--\eqref{i1b} on moving domains proved to exist in \cite{FKNNS} possess also an energy inequality (Theorem \ref{t:ex2}). Using this energy inequality we prove the relative entropy (energy) inequality (Theorem \ref{t:REI}) and finally the weak-strong uniqueness property (Theorem \ref{t:WSU}). %The issue of local existence of strong solutions to the compressible barotropic Navier-Stokes equations on moving domains will be covered by our forthcoming paper \cite{KNP2} and will thus complete the results in this paper.

\section{Preliminaries}\label{p}

For simplicity we assume throughout the paper that $\kappa,\eta = 0$, however all the results hold (with appropriate modifications of formulae) also with $\kappa,\eta > 0$.

We define weak solutions to the compressible Navier-Stokes system on moving domains as follows.
\begin{Definition}\label{d:ws}
We say that the couple $(\vr,\vu)$ is a weak solution of problem \eqref{i1a}-\eqref{i1b} with boundary conditions \eqref{i6}-\eqref{b1} and initial conditions \eqref{i7} if
\begin{itemize}
\item $\vr \in L^\infty(0,T;L^\gamma(\mathbb{R}^3))$, $\vr \geq 0$ a.e. in $Q_T$.
\item $\vu,\nabla_x\vu \in L^2(Q_T)$, $(\vu - \vc{V}) \cdot \vc{n}  (\tau , \cdot)|_{\Gamma_\tau}  = 0$ for a.a. $\tau \in [0,T]$.
\item The continuity equation \eqref{i1a} is satisfied in the whole $\mathbb{R}^3$ provided the density is extended by zero outside $\Omega_t$, i.e.
\bFormula{m1}
\int_{\Omega_\tau} \vr \varphi (\tau, \cdot) \ \dx - \int_{\Omega_0} \vr_0 \varphi (0, \cdot) \ \dx =
\int_0^\tau \int_{ \Omega_t} \left( \vr \partial_t \varphi + \vr \vu \cdot \Grad \varphi \right) \ \dxdt
\eF
for any $\tau \in [0,T]$ and any test function $\varphi \in \DC([0,T] \times \mathbb{R}^3)$.
\item Moreover, equation \eqref{i1a} is also satisfied in the sense of
renormalized solutions introduced by DiPerna and Lions \cite{DL}:
\bFormula{m2}
\int_{\Omega_\tau} b(\vr) \varphi (\tau, \cdot) \ \dx - \int_{\Omega_0} b(\vr_0) \varphi (0, \cdot) \ \dx =
\int_0^\tau \int_{ \Omega_t} \left( b(\vr) \partial_t \varphi + b(\vr) \vu \cdot \Grad \varphi +
\left( b(\vr)  - b'(\vr) \vr \right) \Div \vu \varphi \right) \ \dxdt
\eF
for any $\tau \in [0,T]$, any $\varphi \in \DC([0,T] \times \mathbb{R}^3)$, and any $b \in C^1 [0, \infty)$, $b(0) = 0$, $b'(r) = 0$ for large $r$.
\item The momentum equation \eqref{i1b} is replaced by the family of integral identities
\bFormula{m3}
\int_{\Omega_\tau} \vr \vu \cdot \vph (\tau, \cdot) \ \dx - \int_{\Omega_0} (\vr \vu)_0 \cdot \vph (0, \cdot) \ \dx
\eF
\[
= \int_0^\tau \int_{\Omega_t} \left( \vr \vu \cdot \partial_t \vph + \vr [\vu \otimes \vu] : \Grad \vph + p(\vr) \Div \vph
- \tn{S} (\Grad \vu) : \Grad \vph\right) \dxdt
\]
for any $\tau \in [0,T]$ and any test function $\vph \in \DC([0,T] \times \mathbb{R}^3)$ satisfying
\bFormula{m4}
\vph \cdot \vc{n}|_{\Gamma_\tau} = 0 \ \mbox{for any} \ \tau \in [0,T].
\eF
\end{itemize}
\end{Definition}

The existence of weak solutions to the problem \eqref{i1a}-\eqref{i7} was proved in \cite{FKNNS}.
\begin{Theorem}\label{t:ex}
Let $\Omega_0 \subset \mathbb{R}^3$ be a bounded domain of class $C^{2 + \nu}$, and let $\vc{V} \in C^1([0,T]; C^{3}_c (\mathbb{R}^3))$ be given.
Assume that the pressure $p \in C[0, \infty) \cap C^1(0, \infty)$ satisfies
\[
p(0) = 0,\ p'(\vr) > 0 \ \mbox{for any}\ \vr > 0,\ \lim_{\vr \to \infty} \frac{p'(\vr)}{\vr^{\gamma - 1}} = p_\infty > 0
\ \mbox{for a certain}\ \gamma > 3/2.
\]
Let the initial data satisfy
\[
\vr_0 \in L^\gamma (R^3),\ \vr_0 \geq 0, \ \vr_0 \not\equiv 0,\ \vr_0|_{R^3 \setminus \Omega_0} = 0,\
(\vr \vu)_0 = 0 \ \mbox{a.a. on the set} \ \{ \vr_0 = 0 \} ,\ \int_{\Omega_0} \frac{1}{\vr_0} |(\vr \vu)_0 |^2 \ \dx < \infty.
\]

Then the problem \eqref{i1a}-\eqref{i7} admits a weak solution on any time interval $(0,T)$ in the sense specified through Definition \ref{d:ws}.
\end{Theorem}

\section{Energy inequality}\label{ei}

Our first result is a crucial observation allowing for proving the main theorem about weak-strong uniqueness.

\begin{Theorem}\label{t:ex2}
Let the assumptions of Theorem \ref{t:ex} be satisfied.

Then the problem \eqref{i1a}-\eqref{i7} admits a weak solution on any time interval $(0,T)$ in the sense specified through Definition \ref{d:ws}, satisfying moreover the energy inequality in the following form
\begin{align}\label{eq:EI}
&\int_{\Omega_\tau} \left( \frac{1}{2} \vr |\vu|^2 + H(\vr) \right)(\tau, \cdot) \ \dx +
\frac{1}{2} \int_0^\tau \int_{\Om_t} \mu \left| \Grad \vu + \Grad^t \vu - \frac{2}{3} \Div \vu \tn{I} \right|^2 \ \dxdt \\ \nonumber
\leq &\int_{\Omega_0} \left( \frac{1}{2 \vr_0} |(\vr \vu)_0 |^2 + H(\vr_0) \right) \ \dx + \int_{\Om_\tau} (\vr \vu \cdot \vc{V}) (\tau, \cdot) \ \dx - \int_{\Om_0} (\vr \vu)_0 \cdot \vc{V}(0, \cdot) \ \dx \\ \nonumber
+ &\int_0^\tau \int_{\Om_t} \left( \mu \left(\Grad \vu + \Grad^t \vu - \frac{2}{3} \Div \vu \tn{I} \right) : \Grad \vc{V} - \vr \vu \cdot \partial_t \vc{V} - \vr \vu \otimes \vu : \Grad \vc{V} - p(\vr) \Div \vc{V} \right)  \dxdt.
\end{align}
\end{Theorem}

\bProof
We follow the same series of approximations and penalizations as it is introduced in \cite[Section 3]{FKNNS} in the proof of Theorem \ref{t:ex}. The starting point is thus the modified energy inequality written on the fixed domain $B$ which is a ball large enough such that $\Omega_t \subset B$ for all $t \in [0,T]$ and $\vc{V} = 0$ on $\partial B$, see formula (3.10) in \cite{FKNNS}
\bFormula{p7}
\int_B \left( \frac{1}{2} \vr |\vu|^2 + H(\vr) + \frac{\delta}{\beta - 1} \vr^\beta \right)(\tau, \cdot) \ \dx +
\frac{1}{2} \int_0^\tau \int_B \mu_\omega \left| \Grad \vu + \Grad^t \vu - \frac{2}{3} \Div \vu \tn{I} \right|^2 \ \dxdt
\eF
\[
+
\frac{1}{\ep} \int_0^\tau \int_{\Gamma_t} \left| \left( \vu - \vc{V} \right) \cdot \vc{n} \right|^2 \ {\rm dS}_x \ \dt
\leq \int_B \left( \frac{1}{2 \vr_{0,\delta}} |(\vr \vu)_{0,\delta} |^2 + H(\vr_{0, \delta}) + \frac{\delta}{\beta - 1} \vr_{0, \delta}^\beta \right) \ \dx
\]
\[
+ \int_B \Big(  (\vr \vu \cdot \vc{V}) (\tau, \cdot) - (\vr \vu)_{0,\delta} \cdot \vc{V}(0, \cdot) \Big)
\ \dx
\]
\[
+ \int_0^\tau \int_B \left( \mu_\omega \left(\Grad \vu + \Grad^t \vu - \frac{2}{3} \Div \vu \tn{I} \right) : \Grad \vc{V} - \vr \vu \cdot \partial_t \vc{V} - \vr \vu \otimes \vu : \Grad \vc{V} -
p(\vr) \Div \vc{V} - \frac{ \delta }{\beta - 1} \vr^\beta \Div \vc{V} \right)  \dxdt.
\]

Passing first with $\ep$ to zero, it is not difficult to observe that using the a priori estimates available, all the terms on the right hand side of \eqref{p7} converge to their counterparts. On the left hand side the last term is positive and thus can be omitted. Finally, using the convexity of $\tn{S}(\Grad\vu):\Grad\vu$ we have 
\begin{equation}\label{eq:convex}
	\int_0^T \int_{B} \tn{S}_\omega( \Grad \vu) :\Grad \vu \, \dxdt \leq \liminf\limits_{\ep \to 0}  \int_0^T \int_{B} \tn{S}_\omega(\Grad \vu_\ep) :\Grad \vu_\ep \, \dxdt.
	\end{equation}

Next, passing with $\omega$ to zero, we first observe that all the terms which include the density can be rewritten as the integrals over $\Omega_t$ instead of integrals over $B$ using the fundamental Lemma 4.1 in \cite{FKNNS}. The viscosity term on the right hand side can be treated easily, in particular the integral over $B \setminus \Omega_t$ vanish due to the fact that $\mu_\omega \sil 0$ on this set. On the left hand side we split the integral of the term with stress tensor into two parts, the integral over $B \setminus \Omega_t$ can be omitted since it is positive and on $\Omega_t$ we use the fact that $\mu_\omega = \mu$ is constant and thus we can use again the convexity of $\tn{S}(\Grad\vu):\Grad\vu$ to obtain similar inequality as \eqref{eq:convex}.

Finally, we pass with $\delta$ to zero. Here we use the - nowadays already standard - results of \cite{EF70} to pass to the limit on the right hand side and to adjust the initial conditions, while on the left hand side we use the weak lower semicontinuity of the energy at time $\tau$.

\qed

\section{Relative energy inequality}\label{rei}

Having now already the energy inequality, we can deduce the relative energy inequality in the spirit of \cite{FeJiNo}. Before stating the theorem, we introduce some notation. For a weak solution $(\vr,\vu)$ and a pair of test functions $(r,\vU)$ defined on $Q_T$ we define the relative energy $\Ecal\Big([\vr,\vu]|[r,\vU]\Big)$  as 
\begin{equation}\label{eq:RE}
\Ecal\Big([\vr,\vu]|[r,\vU]\Big)(\tau) = \int_{\Omega_{\tau}} \left( \frac 12 \vr\abs{\vu-\vU}^2 + H(\vr) - H'(r)(\vr-r)-H(r) \right)(\tau,\cdot)\, \dx. 
\end{equation}

We prove the following

\begin{Theorem}\label{t:REI}

Let $(\vr,\vu)$ be a weak solution to the compressible Navier-Stokes system \eqref{i1a}-\eqref{i7} constructed in Theorem \ref{t:ex2}. Then $(\vr,\vu)$ satisfies the following relative energy inequality
\begin{align}\label{eq:REI}
&\Ecal\Big([\vr,\vu]|[r,\vU]\Big)(\tau) + \int_0^\tau\int_{\Om_t} \left(\tn{S}(\Grad \vu) - \tn{S}(\Grad \vU)\right):\left(\Grad \vu - \Grad \vU\right)\,\dxdt \\ \nonumber \leq\, &\Ecal\Big([\vr_0,\vu_0]|[r(0,\cdot),\vU(0,\cdot)]\Big) + \int_0^\tau \Rcal(\vr,\vu,r,\vU)(t) \dt
\end{align}
for a.a. $\tau \in (0,T)$ and any pair of test functions $(r,\vU)$ such that $\vU \in C^{\infty}_c(\overline{Q_T})$, $\vU\cdot\n = \vV\cdot\n$ on $\Gamma_t$ for $t \in [0,T]$, $r \in C^{\infty}_c(\overline{Q_T})$, $r > 0$.
The remainder term $\Rcal$ is given by
\begin{align}\label{eq:R}
\Rcal(\vr,\vu,r,\vU)(t) &= \int_{\Om_t} \vr(\partial_t\vU + \vu\cdot\Grad\vU)\cdot(\vU-\vu) + \tn{S}(\Grad\vU):(\Grad\vU-\Grad\vu) \,\dx \\ \nonumber
&+ \int_{\Omega_t} \Div\vU(p(r)-p(\vr)) + (r-\vr)\partial_tH'(r) + (r\vU-\vr\vu)\cdot\Grad H'(r)\, \dx
\end{align}
\end{Theorem}
\bProof The idea of the proof is the same as in the original paper \cite{FeJiNo}. We will combine the energy inequality \eqref{eq:EI} provided by Theorem \ref{t:ex2} together with weak formulations of the continuity and momentum equations with suitable test functions. Since $\vU$ is not a proper test function in the momentum equation due to its boundary condition, we test the momentum equation with $\vph = \vU-\vV$ to obtain
\begin{align}\label{eq:pr1}
&\int_{\Omega_\tau} \vr \vu \cdot (\vU-\vV) (\tau, \cdot) \ \dx - \int_{\Omega_0} (\vr \vu)_0 \cdot (\vU-\vV) (0, \cdot) \ \dx \\ \nonumber
= &\int_0^\tau \int_{\Omega_t} \left( \vr \vu \cdot \partial_t (\vU-\vV) + \vr [\vu \otimes \vu] : \Grad (\vU-\vV) + p(\vr) \Div (\vU-\vV)
- \tn{S} (\Grad \vu) : \Grad (\vU-\vV)\right)\ \dxdt.
\end{align}
Subtracting \eqref{eq:pr1} from the energy inequality \eqref{eq:EI} we obtain
\begin{align}\label{eq:pr2}
&\int_{\Omega_\tau} \left(\frac 12\vr\abs{\vu}^2 + H(\vr) - \vr \vu \cdot \vU\right) (\tau, \cdot) \ \dx - \int_{\Omega_0} \frac{1}{2 \vr_0} |(\vr \vu)_0 |^2 + H(\vr_0) - (\vr \vu)_0 \cdot \vU (0, \cdot) \ \dx \\ \nonumber
+ &\int_0^\tau \int_{\Omega_t} \tn{S}(\Grad \vu):(\Grad\vu-\Grad\vU) \ \dxdt  \leq \int_0^\tau \int_{\Omega_t} \left( - \vr \vu \cdot \partial_t \vU - \vr [\vu \otimes \vu] : \Grad \vU - p(\vr) \Div \vU \right)\ \dxdt.
\end{align}
Next, we use in the continuity equation as a test function the quantities $\frac 12 \abs{\vU}^2$ and $H'(r)$ respectively to obtain
\bFormula{eq:pr3}
\int_{\Omega_\tau} \frac 12 \vr \abs{\vU}^2 (\tau, \cdot) \ \dx - \int_{\Omega_0} \frac 12 \vr_0 \abs{\vU}^2 (0, \cdot) \ \dx =
\int_0^\tau \int_{ \Omega_t} \left( \vr \vU \cdot \partial_t \vU + \vr \vu \cdot \Grad \vU \cdot \vU \right) \ \dxdt
\eF
and
\bFormula{eq:pr4}
\int_{\Omega_\tau} \vr H'(r) (\tau, \cdot) \ \dx - \int_{\Omega_0} \vr_0 H'(r) (0, \cdot) \ \dx =
\int_0^\tau \int_{ \Omega_t} \left( \vr \partial_t H'(r) + \vr \vu \cdot \Grad H'(r) \right) \ \dxdt
\eF
Adding \eqref{eq:pr3} and subtracting \eqref{eq:pr4} from \eqref{eq:pr2} we obtain
\begin{align}\label{eq:pr5}
&\int_{\Omega_\tau} \left(\frac 12\vr\abs{\vu-\vU}^2 + H(\vr) - H'(r)\vr\right) (\tau, \cdot) \ \dx - \int_{\Omega_0} \frac{1}{2 \vr_0} |(\vr \vu)_0 - \vr_0\vU(0,\cdot) |^2 + H(\vr_0) - H'(r(0,\cdot))\vr_0 \ \dx \\ \nonumber
+ &\int_0^\tau \int_{\Omega_t} \tn{S}(\Grad \vu):(\Grad\vu-\Grad\vU) \ \dxdt  \leq \int_0^\tau \int_{\Omega_t} \left( (\vr \partial_t \vU + \vr\vu\cdot\Grad\vU)\cdot (\vU-\vu) - p(\vr) \Div \vU\right)\ \dxdt \\ \nonumber
- &\int_0^\tau \int_{\Omega_t} \left( \vr \partial_t H'(r) + \vr \vu \cdot \Grad H'(r) \right)\ \dxdt.
\end{align}
Observing that the definition \eqref{p1pb} implies
\begin{equation}\label{eq:Hpid}
p(r) = rH'(r) - H(r),
\end{equation}
we immediately achieve
\begin{equation}\label{eq:Hpid2}
\partial_t p(r) = r\partial_t H'(r).
\end{equation}
Hence, the inequality \eqref{eq:pr5} can be further rewritten as
\begin{align}\label{eq:pr6}
&\int_{\Omega_\tau} \left(\frac 12\vr\abs{\vu-\vU}^2 + H(\vr) - H'(r)\vr\right) (\tau, \cdot) \ \dx - \int_{\Omega_0} \frac{1}{2 \vr_0} |(\vr \vu)_0 - \vr_0\vU(0,\cdot) |^2 + H(\vr_0) - H'(r(0,\cdot))\vr_0 \ \dx \\ \nonumber
+ &\int_0^\tau \int_{\Omega_t} (\tn{S}(\Grad \vu) - \tn{S}(\Grad\vU)):(\Grad\vu-\Grad\vU) \ \dxdt + \int_0^\tau \int_{\Omega_t} \partial_t p(r)\ \dxdt  \\ \nonumber
\leq &\int_0^\tau \int_{\Omega_t} \left( (\vr \partial_t \vU + \vr\vu\cdot\Grad\vU)\cdot (\vU-\vu) + \tn{S}(\Grad\vU):(\Grad\vu-\Grad\vU)\right)\ \dxdt \\ \nonumber
+ &\int_0^\tau \int_{\Omega_t} \left((r - \vr) \partial_t H'(r) - p(\vr) \Div \vU - \vr \vu \cdot \Grad H'(r) \right)\ \dxdt.
\end{align}
Now we claim that the following identity holds
\begin{equation}\label{eq:pr7}
\int_{\Omega_t} p(r)\Div \vU  + r\vU\cdot\Grad H'(r) \ \dx = \int_{\Omega_t} \Div (\vV p(r)) \ \dx.
\end{equation}
Indeed, using the boundary condition $\vU\cdot\n = \vV\cdot\n$ we write
\begin{align}\label{eq:pr8}
\int_{\Omega_t} p(r)\Div\vU \ \dx &= \int_{\Omega_t} p(r)\Div(\vU-\vV)\ \dx + \int_{\Omega_t} p(r)\Div \vV \ \dx \\ \nonumber
&= -\int_{\Omega_t} \vU\cdot\Grad p(r)\ \dx + \int_{\Omega_t}\Div(\vV p(r))\ \dx = -\int_{\Omega_t} r\vU\cdot\Grad H'(r)\ \dx + \int_{\Omega_t}\Div(\vV p(r))\ \dx,
\end{align}
where we used \eqref{eq:Hpid} as well. Adding \eqref{eq:pr7} to \eqref{eq:pr6} we obtain
\begin{align}\label{eq:pr9}
&\int_{\Omega_\tau} \left(\frac 12\vr\abs{\vu-\vU}^2 + H(\vr) - H'(r)\vr\right) (\tau, \cdot) \ \dx - \int_{\Omega_0} \frac{1}{2 \vr_0} |(\vr \vu)_0 - \vr_0\vU(0,\cdot) |^2 + H(\vr_0) - H'(r(0,\cdot))\vr_0 \ \dx \\ \nonumber
+ &\int_0^\tau \int_{\Omega_t} (\tn{S}(\Grad \vu) - \tn{S}(\Grad\vU)):(\Grad\vu-\Grad\vU) \ \dxdt + \int_0^\tau \int_{\Omega_t} (\partial_t p(r) + \Div (\vV p(r)))\ \dxdt  \\ \nonumber
\leq &\int_0^\tau \int_{\Omega_t} \left( (\vr \partial_t \vU + \vr\vu\cdot\Grad\vU)\cdot (\vU-\vu) + \tn{S}(\Grad\vU):(\Grad\vu-\Grad\vU)\right)\ \dxdt \\ \nonumber
+ &\int_0^\tau \int_{\Omega_t} \left((r - \vr) \partial_t H'(r) + (p(r)- p(\vr)) \Div \vU +(r\vU - \vr \vu) \cdot \Grad H'(r) \right)\ \dxdt.
\end{align}
The proof of Theorem \ref{t:REI} is finished observing that standard transport theorem yields the identity
\begin{align}\label{eq:pr10}
&\int_0^\tau \int_{\Omega_t} (\partial_t p(r) + \Div (\vV p(r)))\ \dxdt = \int_0^\tau \frac{{\rm d}}{\dt}\int_{\Omega_t} p(r)\ \dxdt \\ \nonumber
= &\int_{\Omega_\tau} p(r)(\tau,\cdot) \ \dx - \int_{\Omega_0} p(r)(0,\cdot) \ \dx = \int_{\Omega_\tau} (rH'(r)-H(r))(\tau,\cdot) \ \dx - \int_{\Omega_0} (rH'(r)-H(r))(0,\cdot) \ \dx.
\end{align}

\qed

Note that the class of admissible test functions $(r,\vU)$ can be extended by density arguments in a similar manner as in \cite[Section 3.2.2]{FeJiNo}.

\section{Weak-strong uniqueness}\label{wsu}

Finally, we prove the weak-strong uniqueness principle.

\begin{Theorem}\label{t:WSU}
Let $(\vr,\vu)$ be a weak solution to the compressible Navier-Stokes system \eqref{i1a}-\eqref{i7} constructed in Theorem \ref{t:ex2}. Let $(\vrt,\vut)$ be a strong solution to the same problem satisfying
\begin{equation}
0 < \inf_{Q_T} \vrt \leq \sup_{Q_T} \vrt < \infty
\end{equation}
\begin{equation}
\Grad \vrt \in L^2(0,T;L^q(\Omega_t)), \quad \Grad^2 \vut \in L^2(0,T,L^q(\Omega_t))
\end{equation}
with $q > \max\{3;\frac{6\gamma}{5\gamma-6}\}$, and emanating from the same initial data. Then
\begin{equation}
\vr = \vrt, \, \vu = \vut \quad \text{ in } Q_T.
\end{equation}

\end{Theorem}

\bProof The proof follows the same ideas as in \cite{FeJiNo}, however we present it here for completeness of presentation. Plugging in $(r,\vU) = (\vrt,\vut)$ in the relative energy inequality \eqref{eq:REI} we obtain
\begin{align}\label{eq:pr21}
&\Ecal\Big([\vr,\vu]|[\vrt,\vut]\Big)(\tau) + \int_0^\tau\int_{\Om_t} \left(\tn{S}(\Grad \vu) - \tn{S}(\Grad \vut)\right):\left(\Grad \vu - \Grad \vut\right)\,\dxdt \\ \nonumber &\leq \int_0^\tau \int_{\Om_t} \vr(\partial_t\vut + \vu\cdot\Grad\vut)\cdot(\vut-\vu) + \tn{S}(\Grad\vut):(\Grad\vut-\Grad\vu) \,\dxdt \\ \nonumber
&+ \int_0^\tau \int_{\Omega_t} \Div\vut(p(\vrt)-p(\vr)) + (\vrt-\vr)\partial_tH'(\vrt) + (\vrt\vut-\vr\vu)\cdot\Grad H'(\vrt)\, \dx \dt.
\end{align}
Using the strong formulation of the momentum and continuity equations we find out that 
\begin{equation}\label{eq:pr22}
\partial_t\vut + \vut\cdot\Grad\vut = \frac{1}{\vrt} \tn{S}(\Grad\vut) - \Grad H'(\vrt)
\end{equation}
in $Q_T$. Moreover, multiplying the strong formulation of the continuity equation by $H''(\vrt)$ we obtain
\begin{equation}\label{eq:pr23}
\partial_t H'(\vrt) + \vut\cdot\Grad H'(\vrt) = - \Div \vut \vrt H''(\vrt) = - \Div \vut p'(\vrt).
\end{equation}
Finally, integrating by parts we have for a.a. $t \in (0,\tau)$
\begin{equation}\label{eq:pr24}
\int_{\Omega_t} \tn{S}(\Grad\vut) : (\Grad\vut-\Grad\vu)\,\dx = - \int_{\Omega_t} \Div \tn{S}(\Grad\vut) \cdot (\vut-\vu)\,\dx,
\end{equation}
where the boundary integral vanishes due to boundary condition \eqref{b1} and the fact, that $(\vut-\vu)\cdot\n = 0$ on $\Gamma_t$. Combining \eqref{eq:pr21}, \eqref{eq:pr22}, \eqref{eq:pr23} and \eqref{eq:pr24} we arrive to the following version of the relative energy inequality
\begin{align}\label{eq:pr25}
&\Ecal\Big([\vr,\vu]|[\vrt,\vut]\Big)(\tau) + \int_0^\tau\int_{\Om_t} \left(\tn{S}(\Grad \vu) - \tn{S}(\Grad \vut)\right):\left(\Grad \vu - \Grad \vut\right)\,\dxdt \\ \nonumber &\leq \int_0^\tau \int_{\Om_t} \vr(\vu-\vut)\cdot\Grad\vut\cdot(\vut-\vu) - \Div\vut(p(\vr) - p'(\vrt)(\vr-\vrt) - p(\vrt))\,\dxdt \\ \nonumber
&+ \int_0^\tau \int_{\Omega_t} \frac{1}{\vrt}(\vr-\vrt)\Div\tn{S}(\Grad\vut)\cdot(\vut-\vu)\, \dx \dt.
\end{align}
Now we would like to show that all the terms on the right-hand side of \eqref{eq:pr25} can be absorbed by the left hand side and then use the Gronwall lemma. To do that we need the following estimate which can be easily checked
\begin{align}\label{eq:pr26}
H(\vr) - H'(r)(\vr-r) - H(r) &\geq c(r)(\vr-r)^2 \quad \text{ for } \frac{r}{2} < \vr < 2r \\ \nonumber
&\geq c(r)(1+\vr^\gamma) \quad \text{ otherwise}
\end{align}
and also the following Korn-type inequality
\begin{equation}\label{eq:Korn}
\norma{\vz}{W^{1,2}(\Omega_t)} \leq C\norma{\tn{S}(\Grad\vz)}{L^2(\Omega_t)}
\end{equation}
for all $\vz \in W^{1,2}(\Omega_t)$.

Thus, it is not difficult to observe that
\begin{equation}\label{eq:pr27}
\abs{\int_{\Om_t} \vr(\vu-\vut)\cdot\Grad\vut\cdot(\vut-\vu) - \Div\vut(p(\vr) - p'(\vrt)(\vr-\vrt) - p(\vrt))\,\dx} %\\ \nonumber
\leq C\norma{\Grad\vut}{L^\infty(\Omega_t)}\Ecal\Big([\vr,\vu]|[\vrt,\vut]\Big)(t).
\end{equation}
It remains to handle the last term on the right hand side of \eqref{eq:pr25}. We split the integral into three parts, considering first $\vr$ close to $\vrt$, then $\vr$ small and finally $\vr$ large. We have using the H\"older inequality, the Young inequality and \eqref{eq:Korn}
\begin{align}\label{eq:pr28}
&\abs{\int_{\{\vrt/2 \leq \vr \leq 2\vrt\}} \frac{1}{\vrt}(\vr-\vrt)\Div\tn{S}(\Grad\vut)\cdot(\vut-\vu)\,\dx} \\ \nonumber
\leq &C(\delta)\norma{\frac{1}{\vrt}\Div\tn{S}(\Grad\vut)}{L^3(\Omega_t)}^2\int_{\{\vrt/2 \leq \vr \leq 2\vrt\}}(\vr-\vrt)^2\,\dx + \delta\norma{\vut-\vu}{L^6(\Omega_t)}^2 \\ \nonumber
\leq &C(\delta)\norma{\frac{1}{\vrt}\Div\tn{S}(\Grad\vut)}{L^3(\Omega_t)}^2\Ecal\Big([\vr,\vu]|[\vrt,\vut]\Big)(t) + \delta C \norma{\tn{S}(\Grad(\vut)-\Grad(\vu))}{L^2(\Omega_t)}^2.
\end{align}
The last term can be absorbed into the left hand side for $\delta$ small enough, whereas the first term can be treated using the Gronwall lemma.

On the set where $\vr$ is small we can proceed in the following way
\begin{align}\label{eq:pr29}
&\abs{\int_{\{0\leq\vr<\vrt/2\}} \frac{1}{\vrt}(\vr-\vrt)\Div\tn{S}(\Grad\vut)\cdot(\vut-\vu)\,\dx} \\ \nonumber
\leq &\abs{\int_{\{0\leq\vr<\vrt/2\}} \Div\tn{S}(\Grad\vut)\cdot(\vut-\vu)\,\dx} \\ \nonumber
\leq &C(\delta)\norma{\Div\tn{S}(\Grad\vut)}{L^3(\Omega_t)}^2\int_{\{0\leq\vr<\vrt/2\}}1\,\dx + \delta\norma{\vut-\vu}{L^6(\Omega_t)}^2 \\ \nonumber
\leq &C(\delta)\norma{\Div\tn{S}(\Grad\vut)}{L^3(\Omega_t)}^2\Ecal\Big([\vr,\vu]|[\vrt,\vut]\Big)(t) + \delta C \norma{\tn{S}(\Grad(\vut)-\Grad(\vu))}{L^2(\Omega_t)}^2.
\end{align}
Again, the last term is absorbed into the left hand side for $\delta$ small enough and the first term is treated using the Gronwall lemma.

Finally, consider the integral over the set where $\vr$ is large. Here we distinguish two cases. First for $\gamma \leq 2$ we have
\begin{align}\label{eq:pr30}
&\abs{\int_{\{\vr>2\vrt\}} \frac{1}{\vrt}(\vr-\vrt)\Div\tn{S}(\Grad\vut)\cdot(\vut-\vu)\,\dx} \\ \nonumber
\leq &\abs{\int_{\{\vr>2\vrt\}}\vr\frac{\vr-\vrt}{\vr\vrt}\Div\tn{S}(\Grad\vut)\cdot(\vut-\vu)\,\dx} \\ \nonumber
\leq &\abs{\int_{\{\vr>2\vrt\}}\vr\frac{1}{\vrt}\Div\tn{S}(\Grad\vut)\cdot(\vut-\vu)\,\dx} \\ \nonumber
\leq &C(\delta)\norma{\frac{1}{\vrt}\Div\tn{S}(\Grad\vut)}{L^{\frac{6\gamma}{5\gamma-6}}(\Omega_t)}^2\left(\int_{\{\vr>\vrt/2\}}\vr^\gamma\,\dx\right)^{2/\gamma} + \delta\norma{\vut-\vu}{L^6(\Omega_t)}^2 \\ \nonumber
\leq &C(\delta)\norma{\frac{1}{\vrt}\Div\tn{S}(\Grad\vut)}{L^{\frac{6\gamma}{5\gamma-6}}(\Omega_t)}^2\Ecal\Big([\vr,\vu]|[\vrt,\vut]\Big)^{\frac{2}{\gamma}-1}(t)\Ecal\Big([\vr,\vu]|[\vrt,\vut]\Big)(t) + \delta C \norma{\tn{S}(\Grad(\vut)-\Grad(\vu))}{L^2(\Omega_t)}^2.
\end{align}
In this case the power $\frac{2}{\gamma}-1$ is nonnegative and we use also the property $\Ecal\Big([\vr,\vu]|[\vrt,\vut]\Big) \in L^\infty(0,T)$ to proceed further.

For $\gamma > 2$ we have
\begin{align}\label{eq:pr31}
&\abs{\int_{\{\vr>2\vrt\}} \frac{1}{\vrt}(\vr-\vrt)\Div\tn{S}(\Grad\vut)\cdot(\vut-\vu)\,\dx} \\ \nonumber
\leq &\abs{\int_{\{\vr>2\vrt\}}\vr\frac{\vr-\vrt}{\vr\vrt}\Div\tn{S}(\Grad\vut)\cdot(\vut-\vu)\,\dx} \\ \nonumber
\leq &\abs{\int_{\{\vr>2\vrt\}}\vr^{\frac{\gamma}{2}}\frac{1}{\vrt}\Div\tn{S}(\Grad\vut)\cdot(\vut-\vu)\,\dx} \\ \nonumber
\leq &C(\delta)\norma{\frac{1}{\vrt}\Div\tn{S}(\Grad\vut)}{L^3(\Omega_t)}^2\left(\int_{\{\vr>\vrt/2\}}\vr^\gamma\,\dx\right) + \delta\norma{\vut-\vu}{L^6(\Omega_t)}^2 \\ \nonumber
\leq &C(\delta)\norma{\frac{1}{\vrt}\Div\tn{S}(\Grad\vut)}{L^3(\Omega_t)}^2\Ecal\Big([\vr,\vu]|[\vrt,\vut]\Big)(t) + \delta C \norma{\tn{S}(\Grad(\vut)-\Grad(\vu))}{L^2(\Omega_t)}^2.
\end{align}

Altogether we end up with the inequality
\begin{equation}\label{eq:Gron}
\Ecal\Big([\vr,\vu]|[\vrt,\vut]\Big)(\tau) \leq \int_0^\tau h(t)\Ecal\Big([\vr,\vu]|[\vrt,\vut]\Big)(t)\,\dt
\end{equation}
for some $h(t) \in L^1(0,T)$ and the Gronwall lemma finishes the proof.

\qed

\begin{Remark}
The case of nonzero bulk viscosity coefficient $\eta > 0$ in \eqref{i4} as well as the case of nonzero boundary friction coefficient $\kappa > 0$ in \eqref{b1} can be treated by obvious modifications just adding proper integrals to appropriate formulas.
\end{Remark}

\begin{Remark}
The local-in-time existence of solutions belonging to the regularity class of the strong solution in Theorem \ref{t:WSU} will be proved in our forthcoming paper \cite{KNP2}. 
This will show in particular, that the weak-strong uniqueness property stated in Theorem \ref{t:WSU} is in fact a nonempty result.
\end{Remark}

\end{document}